\documentclass[a4paper,10pt]{amsart}
\usepackage{amsmath,amsthm,latexsym,amscd,amssymb}
\usepackage[all]{xy}
\numberwithin{equation}{subsection}

\newtheorem{prop}{Proposition}[section]
\newtheorem{lem}[prop]{Lemma}
\newtheorem{ddd}[prop]{Definition}
\newtheorem{theorem}[prop]{Theorem}

\newcommand{\Gl}{{\rm Gl}} 

\newcommand{\ind}{\mathop{\mbox{\rm ind}}}

\newcommand{\sign}{\mathop{\rm sign}}
\newcommand{\dom}{\mathop{\rm dom}}

\newcommand{\Se}{\mathfrak{AS}}
\newcommand{\SK}{\mathfrak{ASK}}
\newcommand{\SF}{\mathfrak{ASF}}
\newcommand{\F}{\mathfrak{AF}}

\newcommand{\Pj}{{\mathcal P}}

\newcommand{\M}{{\mathcal M}}
\newcommand{\tr}{{\rm tr}}

\newcommand{\Tr}{\mathop{\rm Tr}}

\newcommand{\N}{{\mathcal N}}

\newcommand{\D}{{\mathcal D}}
\newcommand{\U}{{\mathcal U}}

\newcommand{\A}{{\mathcal A}}
\newcommand{\B}{{\mathcal B}}

\newcommand{\C}{C^{\infty}}

\newcommand{\ra}{\partial}

\newcommand{\ten}{\otimes}

\newcommand{\ve}{\varepsilon}

\DeclareMathOperator{\supp}{supp}
\DeclareMathOperator{\spfl}{sf}
\DeclareMathOperator{\Ran}{Ran}
\DeclareMathOperator{\Ker}{Ker}

\def\bbbr{{\rm I\!R}} 
\def\bbbn{{\rm I\!N}} 

\def\bbbc{{\rm I\!C}}

\def\bbbq{{\mathchoice {\setbox0=\hbox{$\displaystyle\rm Q$}\hbox{\raise
0.15\ht0\hbox to0pt{\kern0.4\wd0\vrule height0.8\ht0\hss}\box0}}
{\setbox0=\hbox{$\textstyle\rm Q$}\hbox{\raise
0.15\ht0\hbox to0pt{\kern0.4\wd0\vrule height0.8\ht0\hss}\box0}}
{\setbox0=\hbox{$\scriptstyle\rm Q$}\hbox{\raise
0.15\ht0\hbox to0pt{\kern0.4\wd0\vrule height0.7\ht0\hss}\box0}}
{\setbox0=\hbox{$\scriptscriptstyle\rm Q$}\hbox{\raise
0.15\ht0\hbox to0pt{\kern0.4\wd0\vrule height0.7\ht0\hss}\box0}}}}

\def\bbbz{{\mathchoice {\hbox{$\sf\textstyle Z\kern-0.4em Z$}}
{\hbox{$\sf\textstyle Z\kern-0.4em Z$}}
{\hbox{$\sf\scriptstyle Z\kern-0.3em Z$}}
{\hbox{$\sf\scriptscriptstyle Z\kern-0.2em Z$}}}}

\def\bbbc{{\mathchoice {\setbox0=\hbox{$\displaystyle\rm C$}\hbox{\hbox
to0pt{\kern0.4\wd0\vrule height0.9\ht0\hss}\box0}}
{\setbox0=\hbox{$\textstyle\rm C$}\hbox{\hbox
to0pt{\kern0.4\wd0\vrule height0.9\ht0\hss}\box0}}
{\setbox0=\hbox{$\scriptstyle\rm C$}\hbox{\hbox
to0pt{\kern0.4\wd0\vrule height0.9\ht0\hss}\box0}}
{\setbox0=\hbox{$\scriptscriptstyle\rm C$}\hbox{\hbox
to0pt{\kern0.4\wd0\vrule height0.9\ht0\hss}\box0}}}}

\title{Spectral flow and winding number in von Neumann algebras}
\author{Charlotte Wahl}
\thanks{The author was partially supported by a grant of the Graduiertenkolleg ``Gruppen und Geometrie'', Georg-August-Universit\"at G\"ottingen}
\address{Gottfried Wilhelm Leibniz Bibliothek\\
Nieders\"achsische Landesbibliothek\\
Waterloostr. 8\\
30169 Hannover\\
Germany} 

\email{ac.wahl@web.de}

\begin{document}

\begin{abstract}
We introduce a new topology, weaker than the gap topology, on the space of selfadjoint operators affiliated to a semifinite von Neumann algebra. We define the real-valued spectral flow for a continuous path of selfadjoint Breuer-Fredholm operators in terms of a generalization of the winding number. We compare our definition with Phillips' analytical definition and derive integral formulas for the spectral flow for certain paths of unbounded operators with common domain, generalizing those of Carey--Phillips. Furthermore we prove the homotopy invariance of the real-valued index. As an example we consider invariant symmetric elliptic differential operators on Galois coverings.   
\end{abstract}

\maketitle

\section{Introduction}
Index theory in von Neumann algebras was introduced by Breuer \cite{b1}\cite{b2}. One important application is to the index theory of elliptic invariant differential operators on Galois coverings, whose foundations were laid with Atiyah's $L^2$-index theorem \cite{a}. In the last ten years much research has been devoted to the real-valued spectral flow in a semifinite von Neumann algebra, which generalizes the spectral flow introduced in \cite{aps}. It began with a topological definition of the real-valued spectral flow for loops of bounded selfadjoint Breuer-Fredholm operators due to Perera \cite{pe1}\cite{pe2}. Phillips presented an analytical definition \cite{p}, which works for general continuous paths of bounded selfadjoint Breuer-Fredholm operators and which was used by Phillips and Carey--Phillips to prove integral formulas for the real-valued spectral flow \cite{p}\cite{cp1}\cite{cp2}. We refer to \cite{bcp} for a survey on the analytic approach, on integral formulas, their applications and for further references.

In parallel, new definitions of the (ordinary) spectral flow for paths of unbounded selfadjoint Fredholm operators have been given \cite{blp}\cite{waspflow}. The straightforward way is to define the spectral flow for unbounded operators as the spectral flow of the bounded transform. However, it is often difficult to decide whether or not the bounded transform of a given path depends continuously on the parameter. In \cite{blp} the spectral flow was defined for paths of unbounded selfadjoint Fredholm operators whose resolvents depend continuously on the parameter, or equivalently which are continuous in the gap topology. This definition was further generalized in \cite{waspflow}. In both papers the spectral flow is expressed in terms of a winding number. 

Our first main result in this paper is a new definition of the real-valued spectral flow in terms of a real-valued winding number. The definition applies to paths $(D_t)_{t \in [0,1]}$ of unbounded selfadjoint Breuer-Fredholm operators, for which there is $\ve >0$ such that $\phi(D_t)$ depends continuously on $t$ for all $\phi\in C_c(\bbbr)$ with $\supp \phi \subset [-\ve,\ve]$. Our method is to modify the approach of \cite{waspflow} in order to make it work in a semifinite von Neumann algebra. The second main result is the proof of new integral formulas. In contrast to those of Carey--Phillips mentioned above we do not require the path to be a bounded perturbation of a fixed selfadjoint operator but assume the domain of the operators to be fixed. For example our formulas apply to a path of invariant symmetric elliptic differential operators of positive order on a Galois covering. For the (ordinary) spectral flow related formulas were derived in \cite{waspflfor}. 

The paper is organized as follows:  In \S \ref{top} we define and study a new topology on the space of selfadjoint unbounded operators affiliated to a semifinite von Neumann algebra. In \S \ref{wind} we review the definition of the generalized winding number. In \S \ref{spflow} we use it to define the real-valued spectral flow for a path of selfadjoint unbounded Breuer-Fredholm operators that is continuous in this new topology. In \S \ref{spflind} we express the real-valued index in terms of the real-valued spectral flow. This will imply the homotopy invariance of the real-valued index. In \S \ref{intfor} we derive integral formulas for the spectral flow. In \S \ref{L2ind} we consider paths of invariant symmetric elliptic differential operators on a Galois covering. Here we make use of the theory of elliptic operators over $C^*$-algebras, which was introduced by Mishchenko-Fomenko \cite{mf}.   

In the following a family or path in a topological space need not be continuous. Projections are selfadjoint.  

{\it Acknowledgements:} The author would like to thank Alan Carey for pointing out to her the Laplace transform method on which the proof of Prop. \ref{psum} is based.

\section{A new topology on the space of selfadjoint affiliated operators}
\label{top}

Let $H$ be a Hilbert space. 

We fix a von Neumann algebra $\N \subset B(H)$ endowed with a faithful normal semifinite trace $\tau$. By $K(\N)$ we denote the norm-closed two-sided ideal generated by projections of finite trace and we write $Q(\N)=\N/K(\N)$ and let $\pi:\N \to Q(\N)$ be the projection. 

A closed densely defined operator $D$ on $H$ is affiliated to $\N$ if and only if its bounded transform $F_D=D(1+D^*D)^{-\frac 12}$ is in $\N$. 

Breuer developed in \cite{b1} \cite{b2} an index theory in von Neumann algebras. See \cite[Appendix B]{pr} and \cite[\S 3]{cp} for further generalizations and modifications as needed here.

An operator $D$ is called a Breuer-Fredholm operator (in $\N$) if it is closed, densely defined, affiliated to $\N$ and if $\pi(F_D) \in Q(\N)$ is invertible. The projection onto the kernel of a selfadjoint Breuer-Fredholm operator has finite trace. The index of a Breuer-Fredholm operator $D$ is defined as $$\ind D= \tau(1_{\{0\}}(D^*D))-\tau(1_{\{0\}}(DD^*)) \ .$$
Here $1_{\{0\}}$ denotes the characteristic function of $\{0\} \subset \bbbr$. Furthermore we will write $1_{\ge 0}$ for the characteristic function of the set $\{x\ge 0\} \subset \bbbr$.

We denote the space of selfadjoint operators affiliated to $\N$ by $AS(\N)$ and the subspace of Breuer-Fredholm operators by $ASF(\N)$. We write $AS(\N)_{gap}$ for the set $AS(\N)$ endowed with the gap topology, which is the weakest topology such that the map 
$$AS(\N) \to \N,~D \mapsto (D+i)^{-1}$$ 
is continuous, and we write $ASF(\N)_{gap}$ for the set $ASF(\N)$ with the subspace topology of $AS(\N)_{gap}$. 

The following well-known property of the gap topology is given here for further reference.

\begin{lem}
\label{contgap}
The map
$$C_0(\bbbr)\times AS(\N)_{gap} \to \N,~ (f,D) \mapsto f(D)$$ 
is continuous.

In particular the map 
$$C(\bbbr) \times \N \to \N,~(f,F) \mapsto f(F)$$ 
is continuous. 
\end{lem}

\begin{proof} The first statement follows from $\|f(D)\| \le \sup_{x\in \bbbr} |f(x)|$ and the fact that the algebra generated by the functions $(x \pm i)^{-1}$ is dense in $C_0(\bbbr)$. The second statement follows the continuity of the inclusion $\N \to AS(\N)_{gap}$. 
\end{proof}

Let $\phi \in \C_c(\bbbr)$ be an even function with $\supp \phi = [-1,1]$ and with $\phi'(x)>0$ for $x \in (-1,0)$ and define $\phi_n \in \C_c(\bbbr)$ by $\phi_n(x):=\phi(nx)$ for $n \in \bbbn_0$. For $n \in \bbbn$ let $\psi_n \in C_c(\bbbr)$ be an odd function with values in $[-1,1]$ and with $\psi_n(x)=x$ for $x \in [-\frac 1n, \frac 1n]$ and $\supp \psi_n \subset (-\frac{1}{n-1},\frac{1}{n-1})$ if $n>1$ and $\psi \in C_0(\bbbr)$ if $n=1$. Let $\psi_0(x)=x$.

For $n \in \bbbn_0$ let $\Se_n(\N)$ be the set $AS(\N)$ endowed with the weakest topology such that the maps 
$$\Se_n(\N) \to K(\N),~D \mapsto (\psi_n(D)+i)^{-1}K \ ,$$
$$\Se_n(\N) \to K(\N),~D \mapsto (\psi_n(D)-i)^{-1}K \ ,$$
$$\Se_n(\N) \to \N,~ D \mapsto \phi_n(D)$$
are continuous for all $K \in K(\N)$.

Note that for $n=0$ the last condition is trivial.

Compare this with the definition in \cite{waspflow}: If $\N=B(H)$  with the standard trace, where $H$ is a separable Hilbert space, then the first two conditions are equivalent to the continuity of  $D \mapsto (\psi_n(D)\pm i)^{-1}x$ for all $x \in H$. Hence $\mathfrak{S}_n(H) \to \Se_n(\N)$ is continuous, with $\mathfrak{S}_n(H)$ as defined in \cite{waspflow}. 

\begin{lem}
\label{cont}
For $f \in C_0(\bbbr)$ the map $$\Se_0(\N) \times K(\N) \to K(\N),~ (D,K) \mapsto f(D)K$$ is continuous.
\end{lem}

\begin{proof}
If $f(x)=(x\pm i)^{-1}$ the assertion follows from the definition of $\Se_0(\N)$ and from
\begin{eqnarray*}
\lefteqn{\|f(D)K-f(D')K'\|}\\
&=& \|(f(D)-f(D')K+ f(D')(K-K')\| \\
&\le& \|(f(D)-f(D')K\| + \|f(D')\|\|K-K'\| \ .
\end{eqnarray*} 
The general case follows from the previous one since $K(\N)$ is an ideal and since the algebra generated by the functions $(x+i)^{-1}$ and $(x-i)^{-1}$ is dense in $C_0(\bbbr)$. 
\end{proof}

\begin{lem}
\label{contev}
Let $n\in \bbbn$. For $f \in C_c(\bbbr)$ with $\supp f \subset (-\frac 1n,\frac 1n)$ the map $$\Se_n(\N) \times K(\N) \to K(\N),~ (D,K) \mapsto f(D)K$$ is continuous. If $f$ is even, then the map $$\Se_n(\N) \to \N,~D \mapsto f(D)$$ is continuous.
\end{lem}

\begin{proof}
For any $f \in C_c(\bbbr)$ with $\supp f \subset (-\frac 1n,\frac 1n)$ there is $g \in C_c(\bbbr)$ such that $f=(g\circ \psi_n) \phi_n$. Furthermore $\Se_n(\N) \to \Se_0(\N),~D \mapsto \psi_n(D)$ is continuous. Now the first assertion follows from the previous lemma.

If $f \in C_c(\bbbr)$ is an even function with $\supp f \subset (-\frac 1n,\frac 1n)$, then there is $g \in C_c(\bbbr)$ such that $f=g \circ \phi_n$. Hence by Lemma \ref{contgap} the map $\Se_n \to \N,~ D\mapsto g\circ \phi_n(D)$ is continuous.
\end{proof}

It follows that the identity induces continuous maps $\Se_m(\N) \to \Se_n(\N)$ for $m,n \in \bbbn,~m \le n$. Let $\Se(\N)$ be the set $AS(\N)$ endowed with the direct limit topology. The previous lemma also implies that the definition of $\Se(\N)$ does not depend on the choice of the functions $\phi$ and $\psi_n$.

For $n \in \bbbn$ we define $$ASF_n(\N):=\{D \in ASF(\N)~|~ \phi_n(D) \in K(\N) \}$$ and denote by $\SF_n(\N)$ the set $ASF_n(\N)$ endowed with the subspace topology of $\Se_n(\N)$. 

Furthermore we define $$ASK(\N):=\{D \in AS(\N)~|~ (1+D^2)^{-1} \in K(\N) \} \subset ASF(\N)$$ 
and denote by $\SK(\N)$ the set $ASK(\N)$ endowed with the subspace topology of $AS(\N)_{gap}$. If $D \in ASK(\N)$, then the resolvents of $D$ are in $K(\N)$: From $(1+D^2)^{-1} \in K(\N)$ it follows that $(1+D^2)^{-1/2} \in K(\N)$, thus $(D \pm i)^{-1}\in K(\N)$ since $(D \pm i)^{-1}(1+D^2)^{1/2} \in \N$.

There is a continuous inclusion $\SK(\N) \to \SF_n(\N)$. 

Let $\SF(\N)$ be the inductive limit of the spaces $\SF_n(\N)$. Since $D \in AS(\N)$ is Breuer-Fredholm if and only if $\phi_n(D) \in K(\N)$ for $n$ big enough, the underlying set of $\SF(\N)$ is $ASF(\N)$.

\begin{lem}
\label{contzwo}
\begin{enumerate}
\item For $f\in C_0(\bbbr)$ the map 
$$\SK(\N) \to K(\N),~D \mapsto f(D)$$ 
is continuous.

\item Let $n \in \bbbn$.
For $f\in C_0(\bbbr)$ with $\supp f \subset (-\frac 1n,\frac 1n)$ the map 
$$\SF_n(\N) \to K(\N),~D \mapsto f(D)$$ 
is continuous.
\end{enumerate}
\end{lem}

\begin{proof} (1) Since the algebra generated by the functions $(x\pm i)^{-1}$ is dense in $C_0(\bbbr)$, we have that $f(D)$ in $K(\N)$ if $(D \pm i)^{-1} \in K(\N)$. The continuity follows from Lemma \ref{contgap}.

(2) For $n>0$ there is $g \in C_c(\bbbr)$ with $\supp g \subset (-\frac 1n,\frac 1n)$ such that $f=g \phi_n$. Since $\SF_n(\N) \to K(\N),~D\mapsto \phi_n(D)$ is continuous, the assertion follows from the previous lemma.
\end{proof}

\begin{ddd}
\label{normal}
\begin{enumerate}
\item A {\rm  normalizing} function for $\SK(\N)$ is an odd non-decreasing function $\chi \in C(\bbbr)$ with $\lim_{x \to \infty}\chi (x)=1$ and $\chi^{-1}(0)=\{0\}$.
\item 
A {\rm  normalizing} function for $\SF_n(\N)$ is an odd non-decreasing function $\chi \in C(\bbbr)$ such that $\lim_{x \to \infty}\chi (x)=1$, furthermore $\chi^{-1}(0)=\{0\}$  and $\supp(\chi^2-1) \subset (-\frac 1n,\frac 1n)$.
\end{enumerate}
\end{ddd}

If $\chi$ is a normalizing function for $\SK(\N)$, then $$\SK(\N) \to K(\N),~D \mapsto \chi(D)^2-1$$ is well-defined and continuous by the previous corollary. More generally, for any $f \in C([-1,1])$ with $f(-1)=f(1)=1$ we have that $f \circ \chi-1 \in C_0(\bbbr)$, hence $$\SK(\N) \to K(\N),~D \mapsto f(\chi(D))-1$$ is continuous. The analogous statements for $\SF_n(\N)$ also hold.

Let $B$ be a compact space. If $(D_b)_{b\in B}$ is a continuous family in $\SF(\N)$, then we say that $\chi$ is a normalizing function for $(D_b)_{b \in B}$ if $\chi$ is a normalizing function for $\SF_n(\N)$ with $n$ such that $(D_b)_{b \in B}$ is a continuous family in $\SF_n(\N)$. Here we allow $n=0$ by setting $\SF_0(\N)=\SK(\N)$.

The topological spaces we introduced above have similar properties as the spaces defined in \cite{waspflow}: If  $f: \bbbr \to \bbbr$ is an odd non-decreasing continuous function with $f^{-1}(0)=0$, then $f:\Se(\N) \to \Se(\N)$ can be shown, by using Lemma \ref{contev}, to be continuous. 
   
Furthermore if $(D_b)_{b \in B}$ is a continuous family in $\SF(\N)$ and $(U_b)_{b \in B}$ is a family of unitaries in $\N$ such that $(U_b K)_{b \in B}$ is continuous for each $K \in K(\N)$, then $(U_b D_bU_b^*)_{b \in B}$ is a continuous family in $\SF(\N)$. 

Sometimes we have to vary the function. Then the argument will be similar to the following. If $f_0,f_1$ are odd non-decreasing continuous function with $f_i^{-1}(0)=0,~i=0,1$ and if we set $f_t=(1-t)f_0+tf_1$ for $t \in (0,1)$, then for $D \in ASF(\N)$
$$[0,1] \to \SF(\N),~ t\mapsto f_t(D)$$
is continuous: For $n$ large enough we have that $\psi_n(D) \in K(\N)$ and for $g \in C_c(\bbbr)$ with support small enough we have that $g(f_t(\psi_n(D)))=g(f_t(D))$. By Lemma \ref{contgap} the operator $g(f_t\psi_n(D))$ depends continuously on $t$. Hence for $m$ big $\phi_m(f_t(D))$ and $(\psi_m(f_t(D))\pm i)^{-1}$ depend continuously on $t$.

Compare the following proposition with \cite[Prop. 1.7]{waspflow}. The proof here is analogous. It is given for completeness.
    
\begin{prop}
\begin{enumerate}
\item The identity induces a continuous map from $AS(\N)_{gap}$ to $\Se(\N)$.
\item The set $ASF(\N)$ is open in $\Se(\N)$.
\item The identity induces a homeomorphism from $\Se(\N) \cap ASF(\N)$ to $\SF(\N)$. 
\end{enumerate}
\end{prop}

\begin{proof}
(1) follows from Lemma \ref{contgap}.

(2) Let $D_0 \in ASF_n(\N)$ and let $\chi$ be a normalizing function for $\SF_n(\N)$. Then $\chi(D_0)^2$ is invertible in $Q(\N)$. Since $$\Se_n(\N) \to \N,~D \mapsto \chi(D)^2-1$$ is continuous by Lemma \ref{contev}, the map $\Se_n(\N) \to Q(\N),~D \mapsto  \chi(D)^2$ is continuous as well. Hence there is an open neighbourhood $U$ of $D_0$ in $\Se_n(\N)$ such that $\chi(D)^2$ is invertible in $Q(\N)$ for all $D \in U$. This implies that $U \subset ASF(\N)$. 

(3) is clear.
\end{proof}

The previous constructions generalize to (not necessarily selfadjoint) affiliated operators as follows:

We fix the trace $\Tr \ten \tau$ on the semifinite von Neumann algebra $M_n(\N)=M_n(\bbbc) \ten \N$.

We denote by $A(\N)$ the space of closed densely defined operators affiliated to $\N$, which we consider as a subspace of $AS(M_2(\N))$ via
the map 
$$G:A(\N) \to AS(M_2(\N)),~ D \mapsto \left(\begin{array}{cc}  0 & D^* \\ D & 0 \end{array}\right) \ .$$
Furthermore $AF(\N)$ means the subspace of Breuer-Fredholm operators: We write $AF_n(\N):=AF(\N) \cap ASF_n(M_2(\N))$. Let $\F_n(\N)$ resp. $\F(\N)$ be the space $AF_n(\N)$ resp. $AF(\N)$ with the subspace topology of $\SF_n(M_2(\N))$ resp. $\SF(M_2(\N))$.

We say that $\chi$ is a normalizing function for $\F_n(\N)$ if $\chi$ is a normalizing function for $\SF_n(M_2(\N))$ and then define $\chi(D), \chi(D^*) \in \N$ for $D \in \F_n(\N)$ by $$\chi(G(D))=\left(\begin{array}{cc}  0 & \chi(D^*) \\ \chi(D) & 0 \end{array}\right) \ .$$

{\bf Remark.} In the following we point out one crucial difference between the situation considered here and the special case $\N=B(H)$ for a separable Hilbert space $H$, which was studied in \cite{waspflow}. In the case $\N=B(H)$ we have that $K(\N)=K(H)$. If $(T_b)_{b\in B}$ is a uniformly bounded family in $B(H)$, then the map $$B \to K(H),~b \mapsto T_bK$$ is continuous for all $K \in K(H)$ if and only if $$B\to H,~b \mapsto T_bx$$ is continuous for all $x \in H$. Hence if $B$ is compact and $(D_b)_{b \in B}$ is a continuous family in $\Se_0(B(H))$, then for $f\in C(\bbbr)$ the map
$$B\to K(H),~ b\mapsto f(D_b)K$$ is continuous for all $K \in K(H)$. The last statement does not generalize to an arbitrary semifinite von Neumann algebra $\N$. (By Lemma \ref{cont} it generalizes if we assume that $f \in C_0(\bbbr)$.) The standard proof, see for example \cite[Prop. 1.1]{waspflow}, does not work here since an element $D \in AS(\N)$ is not necessarily densely defined as an operator on the Banach space $K(\N)$.  Differently put, $(D+i)^{-1}K(\N)$ need not be dense in $K(\N)$. This will be demonstrated in the first part of the following example. In the second part we construct a continuous family $(D_b)_{b \in B}$ in $\Se_0(\N)$ and give a $K \in K(\N)$ such that $B \to K(\N),~ b \mapsto D_b(1+D_b^2)^{-1/2}K \in K(\N)$ is not continuous.

The example:

Let $\N=L^{\infty}(\bbbr)$ act on $H=L^2(\bbbr)$ by multiplication. Define the trace $\tau(f)=\int_{-\infty}^{\infty} f(x)~dx$ for $f \in L^{\infty}(\bbbr) \cap L^1(\bbbr)$. Take $g(x)=\frac{1}{|x|}+1$. Clearly $g \in AS(L^{\infty}(\bbbr))$. Since $L^{\infty}(\bbbr) \cap L^1(\bbbr)$ is dense in $K(\N)$, the set $g^{-1}(L^{\infty}(\bbbr) \cap L^1(\bbbr))$ is dense in the domain of $g$ as an operator on $K(\N)$. It is a subset of $L^{\infty}(\bbbr) \cap L^1(\bbbr)$, which is not dense with respect to the $L^{\infty}$-norm, hence $g$ is not densely defined as an operator on $K(\N)$.   

Let $B=\bbbn \cup \{\infty\}$ be the one point compactification of $\bbbn \subset \bbbr$. For $n \in \bbbn$ we define $g_n \in L^{\infty}(\bbbr)$ as follows: We set 
$$g_n(x)=\left\{\begin{array}{l} k \mbox{ if }  |x| > \frac 1n \mbox{ and }|x| \in ]\frac {1}{2k+1}, \frac{1}{2k}] \mbox{ for } k \in \bbbn\\
-k \mbox{ if }|x| > \frac 1n \mbox{ and }|x| \in ]\frac {1}{2k}, \frac{1}{2k-1}] \mbox{ for }k \in \bbbn\\
n \mbox{ for } |x|\in ]0, \frac 1n]\\
0 \mbox{ else.}
\end{array}\right.$$ 
We define $g_{\infty}$ as the pointwise limit of the functions $g_n$. Then the function $(g_n + i)^{-1}$ converges in $L^{\infty}(\bbbr)$ to $(g_{\infty}+i)^{-1}$ for $n \to \infty$, hence the family $(g_n)_{n \in B}$ is continuous in $\Se_0(\N)$, even in $AS(\N)_{gap}$. However for $f(x)=x(1+x^2)^{-1/2}$ and $1_{[-1,1]} \in K(\N)$ the function  $(f\circ g_n)1_{[-1,1]}$ does not converge to $(f\circ g_{\infty})1_{[-1,1]}$ in $L^{\infty}(\bbbr)$ for $n \to \infty$.

\section{Real-valued winding number}
\label{wind}

The definition of the winding number and its properties, which we discuss in this section, are well-known at least in closely related contexts.  The winding number is discussed for Banach algebras with trace in \cite[\S 1]{hs} and appears also in \cite[Appendix \S 4]{c}. We give the proofs for completeness.

We write $l^1(\N)$ for the ideal of operators $A \in \N$ such that $\tau(|A|)<\infty$. Endowed with the norm $\|A\|_1:=\|A\|+\tau(|A|)$ it is a Banach space (see \cite[Prop. A.1]{cp0}) and it holds $$\|SAT\|_1 \le \|S\| \|A\|_1 \|T\|$$ for $S,T \in \N$ and $A \in l^1(\N)$. Moreover $l^1(\N)$ is dense in $K(\N)$.

As mentioned before, the projection onto the kernel of a selfadjoint Breuer-Fredholm operator is in $l^1(\N)$. This implies that any projection $P \in K(\N)$ is in $l^1(\N)$ since the operator $(1-P)$ is Breuer-Fredholm by $\pi(1-P)=\pi(1)\in Q(\N)$. We will need this fact later on.

Let $\Gl(\N) \subset \N$ be the group of invertible elements and let $$\Gl_{K}(\N)=\{1+K \in \Gl(\N)~|~ K \in K(\N) \} \ .$$
Furthermore $\U(\N) \subset \N$ denotes the unitary group and $$\U_{K}(\N):=\U(\N) \cap \Gl_K(\N) \ .$$  
We endow these groups with the subspace topology.   

We define for $s:[a,b] \to \Gl_K(\N)$ with $s-1 \in C^1([a,b],l^1(\N))$
$$w(s):=\frac{1}{2\pi i} \int_a^b \tau(s(x)^{-1}s'(x))~dx \in \bbbr \ .$$ 

For $U \in \Gl_K(\N) \cap (1+l^1(\N))$
$$w(sU)=w(Us)=w(s)$$ and for $U \in \Gl(\N)$ 
$$w(U^{-1}sU)=w(s) \ .$$ 
The following lemma is used to prove homotopy invariance of $w$.

\begin{lem}
If $h:[a,b]\times [c,d] \to \Gl_K(\N)$ such that $h-1 \in C^1([a,b]\times [c,d],l^1(\N))$, then 
$$w(h(\cdot,c))+w(h(b,\cdot))-w(h(a,\cdot))-w(h(\cdot,d))=0 \ .$$ 
\end{lem}

\begin{proof} We can subdivide the rectangle $[a,b]\times [c,d]$ into smaller rectangles and add up the contributions of the pieces without changing the value of the left hand side. Therefore and by compactness it is enough to prove the assertion for $h$ such that in $C^1([a,b]\times [c,d],l^1(\N))$
$$\|h-h(a,a)\|_{C^1} <\|h(a,a)^{-1}\|^{-1} \ .$$ We set $g(x,y)=h(x,y)h(a,a)^{-1}$. It is enough to show the assertion for $g$. 
Since $\|g-1\|_{C^1}<1$ in $C^1([a,b]\times [c,d],l^1(\N))$, the logarithm $f(x,y)=\log(g(x,y))$ is well-defined in $C^1([a,b]\times [c,d],l^1(\N))$. We have that 
\begin{eqnarray*}
w(g(\cdot,c))&=& \frac{1}{2\pi i}  \int_a^b \tau(\ra_x f(x,c))~dx \\
&=& \frac{1}{2\pi i}  \int_a^b \ra_x\tau(f(x,c))~dx \\
&=& \frac{1}{2\pi i} (\tau (f(b,c)) - \tau (f(a,c))) \ .
\end{eqnarray*}
Similar equations hold for the other three terms. Inserting them implies the claim.
\end{proof}

We consider a function on $S^1$ as a function on $[0,1]$ whose endpoints coincide. We denote by $\pi_1(\U_K(\N))$ the free first homotopy group of $\U_K(\N)$.
 
\begin{prop}
The winding number extends to a homomorphism
$$w:\pi_1(\U_{K}(\N)) \to \bbbr \ .$$
\end{prop}

\begin{proof}
Since $C^1(S^1,l^1(\N))$ is dense in $C(S^1,K(\N))$, for any loop $s$ in $\U_K(\N)$ there is a loop $\tilde s$ homotopic to $s$ in $\Gl_K(\N)$ such that $1- \tilde s \in C^1(S^1,l^1(\N))$. We define $w(s):=w(\tilde s)$. This is well-defined: If there is a homotopy in $\Gl_K(\N)$ between two loops $\tilde s_1$ and $\tilde s_2$ with $\tilde s_i-1 \in C^1(S^1,l^1(\N)),~i=1,2$, then there is a homotopy $h$ between $\tilde s_1$ and $\tilde s_2$ in $\Gl_K(\N)$ with $h-1 \in C^1(S^1 \times [0,1],l^1(\N))$. By the previous lemma $w(\tilde s_1)=w(\tilde s_2)$. 
\end{proof}

The definition of the winding number extends to matrices in a straightforward way. The integral formula holds also for $s$ with $(s-1) \in C^1_{pw}(S^1,l^1(\N))$. Here $C^1_{pw}(S^1,l^1(\N)) \subset C(S^1,l^1(\N))$ denotes the space of continuous piecewise differentiable functions with piecewise continuous derivative.

We sketch a more abstract way  via $K$-theory of defining the winding number. 

First we note that $l^1(\N)$ is closed under holomorphic functional calculus in $K(\N)$: If $T \in l^1(\N),~\lambda \in \bbbc$ and $C \in \bbbc+K(\N)$ are such that $(\lambda+T)C=1$, then $C=\lambda^{-1}(1-TC) \in \bbbc+l^1(\N)$.

Since furthermore $l^1(\N)$ is dense in $K(\N)$, we have that $K_*(K(\N)) \cong K_*(l^1(\N))$. Hence the trace induces a homomorphism $\tau:K_0(K(\N)) \to \bbbr$. Furthermore there is a homomorphism 
\begin{eqnarray*}
\pi_1(\U_K(\N)) &\to& K_1(C(S^1,K(\N))) \\
&\cong& K_1(K(\N)) \oplus K_1(C_0((0,1),K(\N))) \\
&\to& K_1(C_0((0,1),K(\N)))\cong K_0(K(\N)) \ ,
\end{eqnarray*}
where the third map is defined by mapping  $K_1(K(\N))$ to zero and the last map is given by Bott periodicity.  It can be checked that the composition of this map with $\tau$ agrees with $w$.

The motivation for introducing the concreter definition of the real-valued winding number in terms of the integral formula is that it can be used to derive new integral formulas for the spectral flow, see \S \ref{intfor}. (See also \cite{waspflfor} for the case $\N=B(H)$ for a separable Hilbert space $H$.)

\section{Real-valued spectral flow}
\label{spflow}

Before defining the real-valued spectral flow in terms of the winding number we review Phillips' analytic definition of the real-valued spectral flow in a semifinite von Neumann algebra, which we call analytic spectral flow, denoted by $\spfl^a$, for the moment. It applies to any path $(D_t)_{t \in [0,1]}$ in $ASF(\N)$ such that the bounded transform $F_{D_t}$ depends continuously on $t$. Let $P_t:=1_{\ge 0}(D_t)$. Recall the projection $\pi:\N \to Q(\N)$. If $\|\pi(P_s)-\pi(P_{t})\| < \frac 12$ for all $s,t \in [0,1]$, then $$\spfl^a((D_t)_{t \in [0,1]}):=\ind(P_0P_1) \ .$$ This is well-defined since $P_0P_1$ is Breuer-Fredholm in $P_0\N P_1$. The general case can be reduced to this case by cutting the path into small enough pieces and adding up the contributions of the pieces.  This works since $\pi(P_t)$ depends continuously on $t$: Let $\chi$ be a normalizing function of the path $(D_t)_{t\in [0,1]}$. Then $\pi(P_t)=\pi(\frac 12(\chi(D_t)+1))$. Continuity of the bounded transform implies that $\chi(D_t)$ depends continuously on $t$, hence also $\pi(P_t)$.

\begin{ddd}
Let $(D_t)_{t\in [0,1]}$ be a continuous path in $\SF(\N)$ with invertible endpoints. Let $\chi \in C(\bbbr)$ be a normalizing function of $(D_t)_{t \in [0,1]}$  such that $\chi(D_0)$ and $\chi(D_1)$ are involutions. Then we define
$$\spfl((D_t)_{t \in [0,1]}):=w((e^{\pi i (\chi(D_t)+1)})_{t \in [0,1]}) \in \bbbr \ .$$ 
\end{ddd}

The right hand side is well defined since $e^{\pi i (\chi(D_t)+1)} \in C(S^1,\U_K(\N))$ by the remark after Def. \ref{normal}. 
The homotopy invariance of the winding number implies that the definition does not depend on the choice of the normalizing function: For two normalizing functions $\chi_0,\chi_1$ the interpolation $\chi_s=s\chi_1+(1-s)\chi_0$ is also a normalizing function and $e^{\pi i (\chi_s(D_t)+1)}$ depends continuous on $s$ and $t$.
 
In a similar way the spectral flow can be defined in a $C^*$-algebraic setting as studied in \cite{mu}. Instead of operators affiliated to a von Neumann algebra one would have to consider regular unbounded multipliers on the ideal that plays the role of $K(\N)$. Note that elements of $AS(\N)$ may fail to be regular as unbounded multipliers on $K(\N)$ since they need not be densely defined (see the Remark at the end of \S \ref{top}).

Some properties are in order:

\begin{itemize}
\item If each $D_t$ has a bounded inverse, then $\spfl((D_t)_{t\in [0,1]})=0 \ :$ 

There is $\ve>0$ such that $(-\ve,\ve)$ is a subset of the resolvent set of $D_t$ for all $t \in [0,1]$. Hence for a normalizing function $\chi$ of $(D_t)_{t\in [0,1]}$ with $\supp (\chi^2-1) \subset (-\ve,\ve)$ we have that $e^{\pi i(\chi(D_t)+1)}=1$. Thus $$w((e^{\pi i(\chi(D_t)+1)})_{t \in [0,1]})=0 \ .$$  

\item Let $(U_t)_{t \in [0,1]}$ be a path of unitaries in $\N$ such that $[0,1] \to K(\N),~t \mapsto U_tK$ is continuous for each $K \in K(\N)$. Then $$\spfl((U_tD_tU_t^*)_{t \in [0,1]})= \spfl((D_t)_{t \in [0,1]})\ :$$
Since $$[0,1]\times S^1 \to \U_K(\N),~(s,t) \mapsto U_{st}e^{\pi i (\chi(D_t)+1)}U_{st}^*$$ is continuous, homotopy invariance of the winding number implies that
\begin{eqnarray*}
w((U_te^{\pi i (\chi(D_t)+1)}U_{t}^*)_{t \in [0,1]}) &=& w((U_0e^{\pi i (\chi(D_t)+1)}U_0^*)_{t \in [0,1]})\\
&=&w((e^{\pi i (\chi(D_t)+1)})_{t\in [0,1]}) \ .
\end{eqnarray*} 

\item The spectral flow is homotopy invariant: If $(D_{st})_{(s,t) \in [0,1] \times [0,1]}$ is a continuous family in $\SF(\N)$ and $D_{s0}, D_{s1}$ are invertible for each $s \in [0,1]$, then 
$$\spfl((D_{0t})_{t \in [0,1]})=\spfl((D_{1t})_{t \in [0,1]}) \ .$$
\end{itemize}

The definition of the real-valued spectral flow generalizes to paths with not necessarily invertible endpoints as follows:

If for example $D_0$ is not invertible, then let $\ve \in (0,1)$ be such that $Q=1_{[-\ve,\ve]}(D_0) \in K(\N)$. Since $D_0+Q$ is invertible, we can define $$\spfl((D_t)_{t \in [0,1]}):=\spfl((D_t+(1-t)Q_0)_{t \in [0,1]})-  \spfl^a((D_0+(1-t)Q_0)_{t \in [0,1]})  \ .$$
The definition in the case where $D_1$ is not invertible resp. $D_0,D_1$ are both not invertible is similar.  

By using reparametrisation the spectral flow can be defined also for compact intervals different from $[0,1]$.
It is additive with respect to concatenation of paths.
  
The following technical lemma is needed in the proof of the subsequent proposition. 

\begin{lem}
If $F \in \N$ is a bounded selfadjoint operator such that $\pi(F^2)=1 \in Q(\N)$, then $F+A$ is a Breuer-Fredholm operator for any selfadjoint $A \in \N$ with $\|\pi(A)\|<1$.

Furthermore for $\|\pi(A)\| < \frac 12$ $$\|\pi(1_{\ge 0}(F+A))-\pi(1_{\ge 0}(F))\| <2\|\pi(A)\| \ .$$
\end{lem}

\begin{proof}
Since the spectrum of $\pi(F)$ is a subset of $\{-1,1\}$ and $\|\pi(A)\|<1$, the spectrum of $\pi(F+A)$ does not contain zero. This proves the first assertion.

Now let $\Gamma=\{|\lambda -1| =1\} \subset \bbbc$ with counterclockwise orientation. Then
\begin{eqnarray*}
\lefteqn{\|\pi(1_{\ge 0}(F+A))-\pi(1_{\ge 0}(F))\|}\\
&\le& \|\frac{1}{2\pi i}\int_{\Gamma}(\pi(F+A)-\lambda)^{-1}-(\pi(F)-\lambda)^{-1}~d\lambda \| \\
&=& \|\frac{-1}{2\pi i}\int_{\Gamma}(\pi(F+A)-\lambda)^{-1}\pi(A)(\pi(F)-\lambda)^{-1}~d\lambda \|\\
&\le&\frac{1}{2\pi} \|\pi(A)\| \int_{\Gamma}\|(\pi(F+A)-\lambda)^{-1}\|~ \|(\pi(F)-\lambda)^{-1}\|~ d\lambda\\
&< & 2\|\pi(A) \| \ .
\end{eqnarray*}
\end{proof}

\begin{prop}
\label{eqanaltop}

Let $(D_t)_{t \in [0,1]}$ be a path in $ASF(\N)$ such that $(F_{D_t})_{t \in [0,1]}$ is a continuous path in $\N$. 

Then 
$$\spfl^a((D_t)_{t \in [0,1]})=\spfl((D_t)_{t \in [0,1]}) \ .$$ 
\end{prop}

\begin{proof}
Without loss of generality (by cutting the path into small enough pieces) we may assume that $\|\pi(P_0)-\pi(P_t)\| < \frac {1}{16}$ for all $t \in [0,1]$ where $P_t=1_{\ge 0}(D_t)$. Furthermore it is enough to consider the case where $D_0,D_1$ are invertible. Define $B_t:=(1-t)(2P_0-1)+t(2P_1-1)$. Then $P_i=1_{\ge 0}(B_i),~i=0,1$.
Since $$\|\pi(B_t-B_0)\| = \|\pi(2t(P_1-P_0))\|<\frac 18 \ ,$$ 
by the lemma $$\|\pi(1_{\ge 0}(B_t))- \pi(P_0)\| < \frac 14 \ .$$
Thus for all $s,t \in [0,1]$
\begin{eqnarray*}
\lefteqn{\|\pi(1_{\ge 0}(B_s))- \pi(1_{\ge 0}(B_t))\|}\\
&\le& \|\pi(1_{\ge 0}(B_s))- \pi(P_0)\|~ + ~\| \pi(P_0)- \pi(1_{\ge 0}(B_t))\|  \\ 
&<& \frac 12 \ . 
\end{eqnarray*}

By the definition of the analytic spectral flow 
$$\spfl^a((D_t)_{t\in [0,1]}) = \ind(P_0P_1) = \spfl^a((B_t)_{t \in [0,1]}) \ .$$

Let $\chi$ be a normalizing function of $(D_t)_{t \in [0,1]}$. The map $$[0,1]\times [0,1] \to \SF(\N),~(s,t) \to (1-s)D_t + s\chi(D_t)$$ is continuous, hence by homotopy invariance $$\spfl((D_t)_{t\in [0,1]})=\spfl((\chi(D_t)_{t \in [0,1]}) \ .$$ By the previous lemma the map 
$$[0,1]\times [0,1] \to \SF(\N),~(s,t) \to (1-s)\chi(D_t) + sB_t$$ is well-defined since
$$(1-s)\chi(D_t) + sB_t=\chi(D_t)+s(B_t-\chi(D_t))$$ and
$$\|\pi(B_t-\chi(D_t))\|\le \|\pi(B_t-B_0)\|+ \|\pi(2P_0-1)-\pi(2P_t-1)\| < 1 \ .$$
Furthermore it is continuous.

Again by homotopy invariance this implies  that $$\spfl((\chi(D_t))_{t \in [0,1]})=\spfl((B_t)_{t \in [0,1]}) \ .$$

We will use an idea of \cite[\S 5]{bcp}, where the spectral flow of a path between involutions is discussed, in order to show that $$\spfl((B_t)_{t \in [0,1]})=\spfl^a((B_t)_{t \in [0,1]}) \ .$$ 
Then the assertion follows. 

In the following the (analytic) spectral flow will sometimes be taken in subalgebras of $\N$ without change of notation. The operator $B_t$ is invertible for $t \neq \frac 12$. Choose $r>0$ such that $P:=1_{[0,r]}(B_{1/2}^2) \in l^1(\N)$. This is possible since $B_{1/2}^2$ is a selfadjoint Breuer-Fredholm operator. Since $B_{1/2}^2B_t=B_tB_{1/2}^2$ for all $t \in [0,1]$, we have that $[P,B_t]=0$. Hence $B_t=(1-P)B_t(1-P)+PB_tP$. 
It is clear that $\spfl^a$ and $\spfl$ are additive with respect to direct sums. Since $(1-P)B_t(1-P)$ is invertible for any $t \in [0,1]$, the contribution of $(1-P)\N (1-P)$ to $\spfl^a((B_t)_{t \in [0,1]})$ resp. $\spfl((B_t)_{t\in [0,1]})$ vanishes. We calculate the contribution of the finite von Neumann algebra $P\N P$: By \cite[\S 5.1]{bcp} $$\spfl^a((PB_tP)_{t\in [0,1]})=\frac 12\tau(PB_1P -PB_0P)= \frac 12\int_0^1 \tau(\frac{d}{dt}PB_tP) \ .$$
By the integral formula for the winding number 
\begin{eqnarray*}
\spfl((PB_tP)_{t \in [0,1]})&=&\frac{1}{2\pi i}\int_0^1 \tau(e^{-\pi i(PB_tP+1)}\frac{d}{dt}e^{\pi i(PB_tP+1)})~dt\\
&=& \frac 12\int_0^1 \tau(\frac{d}{dt}PB_tP)~dt
\end{eqnarray*}
The last equality can be proven for example by slightly modifying the proof of \cite[\S I.6.11 Lemma 3]{d}, which uses Cauchy formula. It also follows fromr the results in \S \ref{intfor}. 

This shows the assertion.

An alternative proof of the last part would be to consider the limit $r \to 0$; compare with the proof of Prop. \ref{lemmaspflind}.
\end{proof}

\section{Real-valued spectral flow and index}
\label{spflind}

The real-valued spectral flow and the real-valued index are related as follows:

\begin{prop}
\label{lemmaspflind}
Let $D \in AF(\N)$. Then $$\spfl(\left(\begin{array}{cc} t-\frac 12 & D^* \\ D & \frac 12-t \end{array}\right)_{t \in [0,1]})=\ind(D) \ .$$
\end{prop}

\begin{proof}
Let $$\tilde D_t:=\left(\begin{array}{cc} t-\frac 12 & D^* \\ D & \frac 12-t \end{array}\right) \ .$$
From $$\tilde D_t^2=\left(\begin{array}{cc} (t-\frac 12)^2+D^*D & 0 \\ 0 & (t-\frac 12)^2+DD^* \end{array}\right)$$ it follows that $\tilde D_t$ is invertible for $t \neq \frac 12$. 
We define the Breuer-Fredholm operator $$B_t:=\tilde D_t |\tilde D_0|^{-1} \ .$$ Since $\tilde D_0^2=\tilde D_1^2$, the operators $B_0$ and $B_1$ are involutions. By $\tilde D_t=(1-t)\tilde D_0+ t \tilde D_1$ we have that $$B_t=(1-t)B_0 +t B_1 \ .$$ 
Furthermore $(s,t) \mapsto s\tilde D_t(1+ \tilde D_t^2)^{-\frac 12} + (1-s)B_t$ is a homotopy in $\SF(\N)$ between $(B_t)_{t \in [0,1]}$ and the bounded transform of $(\tilde D_t)_{t \in [0,1]}$. Thus 
$$\spfl((\tilde D_t)_{t \in [0,1]})= \spfl((B_t)_{t \in [0,1]}) \ .$$  
The following argumentation is essentially as in \cite[proof of Prop. 5.3]{bcp}; compare also with the proof of Prop. \ref{eqanaltop}.

There is $r_0>0$ such that for all $0\le r<r_0$ 
$$P_r:=1_{[0,r]}(B_{1/2}^2) \in l^1(\N) \ .$$ 
Furthermore $[P_r,B_t]=0$. Since $P_r$ converges strongly to $P_0$ for $r \to 0$, the trace $\tau(P_r)$ converges to $\tau(P_0)$. Let $r \in (0,r_0)$. The spectral flow of $(B_t)_{t\in [0,1]}$ in $\N$ equals the sum of the spectral flow of $((1-P_r)B_t(1-P_r))_{t \in [0,1]}$ in $(1-P_r)\N(1-P_r)$, the spectral flow of $((P_r-P_0)B_t(P_r-P_0))_{t\in [0,1]}$ in $(P_r-P_0)\N(P_r-P_0)$ and the spectral flow of $(P_0B_tP_0)_{t \in [0,1]}$ in $P_0\N P_0$. The contribution of $(1-P_r)\N(1-P_r)$ vanishes since $(1-P_r)B_t(1-P_r)$ is invertible for all $t \in [0,1]$. 

From $$\spfl((P_r-P_0)B_t(P_r-P_0))_{t \in [0,1]})=\frac 12\tau((P_r-P_0)(B_1-B_0)(P_r-P_0))$$ (see the formula at the end of the proof of Prop. \ref{eqanaltop}), it follows that
$$|\spfl((P_r-P_0)B_t(P_r-P_0))_{t \in [0,1]})| \le \tau(P_r-P_0) \ .$$
Hence the contribution of $(P_r-P_0)\N(P_r-P_0)$ converges to zero for $r \to \infty$.  

Furthermore $$\Ran P_0 =\Ker \tilde D_{1/2}^2=\Ker D \oplus \Ker D^*$$ and with respect to this decomposition
$$P_0B_tP_0=(t-\tfrac 12) \oplus (\tfrac 12-t) \ .$$ Hence $\spfl((P_0B_tP_0)_{t \in [0,1]})=\ind(D)$.
\end{proof}

\begin{prop}
The index $$\ind:\F(\N) \to \bbbr$$ is homotopy invariant.
\end{prop}

\begin{proof}
If $(D_t)_{t \in [0,1]}$ is a continuous path in $\F(\N)$, then
$$[0,1] \times [0,1]  \to \SF(M_2(\N)),~(s,t) \mapsto \left(\begin{array}{cc} t-\frac 12 & D_s^* \\ D_s & \frac 12-t \end{array}\right)$$ is continuous. Now the assertion follows from the previous proposition and the homotopy invariance of the spectral flow.
\end{proof}

\section{Integral formulas}
\label{intfor}

In this section we derive integral formulas for the spectral flow in the spirit of those proven in \cite{cp1}\cite{cp2}. For simplicity we restrict to paths with invertible endpoints. Formulas for paths whose endpoints are not both invertible can then be obtained by using the definition of the spectral flow in \S \ref{spflow}.

We begin by studying how a path of operators $(F_t)_{t\in [0,1]}$ in $\N$ with $\|F_t\| \le 1$ and such that $(t \mapsto F_tK) \in  C^1([0,1],l^1(\N))$ for $K \in l^1(\N)$ behaves under functional calculus. 

Since, as an operator on $l^1(\N)$, the derivative $\frac{d}{dt}F_t$ is the strong limit of the bounded operators $\frac 1h(F_{t+h}-F_t)$ for $h \to 0$, it is a bounded operator on $l^1(\N)$.

\begin{lem}
\label{hol}
 Let $(F_t)_{t \in [0,1]}$ be a path of selfadjoint operators in $\N$ with $\|F_t\| \le 1$ for all $t \in [0,1]$ and such that $(t \mapsto F_tK) \in  C^1([0,1],l^1(\N))$ for all $K \in l^1(\N)$. 

Let $g$ be a holomorphic function defined on a neighbourhood of $[-1,1]$. 
Then for all $K \in l^1(\N)$ 
$$(t \mapsto g(F_t)K) \in C^1([0,1],l^1(\N))$$ and
$$\frac{d}{dt}g(F_t)K= -\int_{\Gamma}g(\lambda) (F_t-\lambda)^{-1}(\frac{d}{dt}F_t)(F_t-\lambda)^{-1}K~ d\lambda \ ,$$
where $\Gamma$ is a closed curve in the domain of $g$, not intersecting $[-1,1]$ and with winding number $1$ with respect to the origin.
\end{lem}

\begin{proof}
By $$((F_{t+h}-\lambda)^{-1}-(F_t - \mu)^{-1})K= - (F_{t+h}-\lambda)^{-1}(F_{t+h}-F_t-\lambda+\mu)(F_t-\mu)^{-1}K \ ,$$
the function $$(\bbbc \setminus [-1,1]) \times [0,1] \to l^1(\N),~(\lambda,t) \mapsto (F_t-\lambda)^{-1}K$$
is continuous. Hence Cauchy formula implies that $(t \mapsto g(F_t)K) \in  C([0,1],l^1(\N))$.
The assertion on the derivative follows from Cauchy formula and $$\frac{d}{dt}(F_t-\lambda)^{-1}K=(F_t-\lambda)^{-1}(\frac{d}{dt}F_t)(F_t-\lambda)^{-1}K \ .$$
\end{proof}

\begin{prop}
\label{der}
Let $(F_t)_{t \in [0,1]}$ be a path of selfadjoint operators in $\N$ with $\|F_t\| \le 1$ for all $t \in [0,1]$ and such that $(t \mapsto F_tK) \in  C^1([0,1],l^1(\N))$ for all $K \in l^1(\N)$. 
Let $g \in C^2_c(\bbbr)$.
 
Then for all $K \in l^1(\N)$ 
$$(t \mapsto g(F_t)K) \in C^1([0,1],l^1(\N))$$ and
$$\frac{d}{dt}g(F_t)K=\int_{-\infty}^{\infty} \int_0^1 i \lambda \hat g(\lambda) e^{ i (1-u)\lambda F_t}(\frac{d}{dt}F_t)e^{ i u \lambda F_t}K~ du~d\lambda \ ,$$
where $\hat g \in L^1(\bbbr)$ is defined via Fourier transform such that
$$g(x)=\int_{-\infty}^{\infty} \hat g(\lambda)e^{ i\lambda x}~d\lambda \ .$$
The integral converges in $l^1(\N)$.
\end{prop}

\begin{proof}
Since the Fourier transform maps $C_c^1(\bbbr)$ to $L^1(\bbbr)$, it holds that $\lambda \hat g(\lambda) \in L^1(\bbbr)$.

We have that 
$$g(F_t)=\int_{-\infty}^{\infty} \hat g(\lambda)e^{ i \lambda F_t} ~d\lambda \ .$$
By the previous lemma $$(t \mapsto e^{ i \lambda F_t}K) \in C^1([0,1],l^1(\N))$$ for any $K \in l^1(\N)$. The expression depends also continuously on $\lambda$. Hence $(t \mapsto g(F_t)K) \in C([0,1],l^1(\N))$.
 
Let $G_s:=sF_{t +h}+(1-s)F_t$. We have that 
\begin{eqnarray*}
e^{ i \lambda F_{t+h}}-e^{i \lambda F_t}&=&\int_0^1 \frac{d}{ds} e^{i \lambda G_s}~ds\\
&=&i\lambda \int_0^1\int_0^1 e^{i(1-u)\lambda G_s}(F_{t+h}-F_t)e^{i u \lambda G_s} \ ~ds ~du \ .
\end{eqnarray*}

Thus for $K \in l^1(\N)$
\begin{eqnarray*}
\frac{d}{dt}e^{i \lambda F_t}K &=&\lim\limits_{h \to 0}\frac{1}{h}(e^{i \lambda F_{t+h}}-e^{i \lambda F_t})K\\
&=&\lim\limits_{h \to 0}\frac{i \lambda}{h} \int_0^1\int_0^1 e^{i(1-u)\lambda G_s}(F_{t+h}-F_t)e^{i u \lambda G_s}K~ds ~du \\
&=&i\lambda \int_0^1 e^{i(1-u)\lambda F_t} (\frac{d}{dt} F_t)e^{i u \lambda F_t}K~du \ .
\end{eqnarray*}
\end{proof}

\begin{lem}
\label{holtrace}
Let $(F_t)_{t \in [0,1]}$ be a path of selfadjoint operators in $\N$ with $\|F_t\| \le 1$ for all $t \in [0,1]$ and such that $(t \mapsto F_tK) \in  C^1([0,1],l^1(\N))$ for all $K \in l^1(\N)$, furthermore such that $(t\mapsto (1-F_t^2)) \in C^1([0,1],l^1(\N))$. If $f$ is a holomorphic function defined on a neighbourhood of $[-1,1]$ with $f(-1)=f(1)=0$, then $(t \mapsto f(F_t)) \in C^1([0,1],l^1(\N))$.
\end{lem}

\begin{proof} We apply Lemma \ref{hol}.
The function $g(z)=f(z)/(z^2-1)$ is holomorphic on the domain of $f$. Hence $f(F_t)=g(F_t)(F_t^2-1) \in C([0,1],l^1(\N))$ and $$\frac{d}{dt}f(F_t)=(\frac{d}{dt}g(F_t))(F_t^2-1) + g(F_t)\frac{d}{dt}(F_t^2-1) \in C([0,1],l^1(\N)) \ .$$
\end{proof}

In the following let $\chi \in C^2_c(\bbbr)$ be an odd function such that $\chi(1)=1$ and $\chi'(0)>0$ and such that $\chi|_{[-1,1]}$ is non-decreasing. Note that $\chi(\frac{x}{(x^2+1)^{1/2}})$ is a normalizing function for $\SK(\N)$.

In the proof of the theorem we will use that if $D \in ASK(\N)$, then $\phi(D) \in l^1(\N)$ for all $\phi \in C_c(\bbbr)$. This can be seen as follows: Let $n \in \bbbn$ be such that $\supp \phi \subset [-n,n]$ and let $\psi \in C_c(\bbbr)$ be such that $\psi|_{[-n,n]}=1$. Since $\psi(D) \in K(\N)$ (this is a special case of Lemma \ref{contzwo}), also $1_{[-n,n]}(D)= 1_{[-n,n]}(D)\psi(D) \in K(\N)$, thus $1_{[-n,n]}(D) \in l^1(\N)$. Hence  $\phi(D)=1_{[-n,n]}(D)\phi(D) \in l^1(\N)$. 
 
\begin{theorem}
\label{thintfor1}
Let $(F_t)_{t \in [0,1]}$ be a path of selfadjoint Breuer-Fredholm operators in $\N$ such that $F_t$ is the bounded transform of an element of $ASK(\N)$ for each $t \in [0,1]$. Assume that $(t \mapsto F_tK) \in  C^1([0,1],l^1(\N))$ for all $K \in l^1(\N)$, furthermore that $(t\mapsto \chi'(F_t)) \in C([0,1],l^1(\N))$ and that $(t\mapsto (\chi(F_t)^2-1)) \in C^1([0,1],l^1(\N))$. Let $F_0,F_1$ be invertible.  

Then 
\begin{eqnarray*}
\spfl((F_t)_{t \in [0,1]})&=&\frac 12 \int_0^1  \tau ((\frac{d}{dt}F_t)\chi'(F_t))~ dt\\
&& +~ \frac 12 \tau(2P_1-1- \chi(F_1)) - \frac 12 \tau(2P_0-1-\chi(F_0))\ ,
\end{eqnarray*} 
where $P_i= 1_{\ge 0}(F_i)$.
\end{theorem}

\begin{proof} 
Since $(2P_i-1)-\chi(F_i)=((2P_i-1)+\chi(F_i))^{-1}(1-\chi(F_i)^2) \in l^1(\N)$, the last two terms on the right hand side are well-defined.

The operators $G_t:=(1-t)(2P_0-1) + t\chi(F_0)$ and $H_t:=(1-t)\chi(F_1)+ t(2P_1-1)$ are invertible for each $t \in [0,1]$, hence  
$$\spfl((G_t)_{t \in [0,1]})=\spfl((H_t)_{t \in [0,1]})=0 \ .$$ 
From 
\begin{eqnarray*}
G_t^2-1 &=& \bigl((2P_0-1) + t(\chi(F_0)-(2P_0-1))\bigr)^2-1\\
&=& t^2(\chi(F_0)-(2P_0-1))^2+2t(2P_0-1)(\chi(F_0)-(2P_0-1))
\end{eqnarray*}
we see that $(t \mapsto (G_t^2-1)) \in C^1([0,1],l^1(\N))$. Analogously $(t \mapsto (H_t^2-1)) \in C^1([0,1],l^1(\N))$.
 
Define $(Q_t)_{t \in [-1,2]}$ by $Q_t:= \chi(F_t)$ for $t \in [0,1]$ and $Q_t:=G_{t+1}$ for $t \in [-1,0]$ and $Q_t:=H_{t-1}$ for $t \in [1,2]$. Then $Q_{-1},Q_2$ are involutions.

By the previous lemma $(t \mapsto (e^{\pi i (Q_t+1)}-1)) \in C_{pw}^1([-1,2],l^1(\N))$.

Set
$$L=\frac{1}{2 \pi i} \int_0^1 \tau  (e^{-\pi i(G_t+1)}\frac{d}{dt}e^{\pi i (G_t+1)})~dt + \frac{1}{2 \pi i} \int_0^1 \tau  (e^{-\pi i(H_t+1)}\frac{d}{dt}e^{\pi i (H_t+1)})~dt\ .$$

The definition of the spectral flow and the integral formula for the winding number imply that
\begin{eqnarray*}
\spfl((F_t)_{t \in [0,1]})&=&
\spfl((Q_t)_{t \in [-1,2]}) \\ 
&=&\frac{1}{2 \pi i} \int_0^1 \tau  (e^{-\pi i(\chi(F_t)+1)}\frac{d}{dt}e^{\pi i ( \chi(F_t)+1)})~dt + L \ .
\end{eqnarray*}

We calculate the first term in the last line.

Let $(\phi_n)_{n \in \bbbn}$ be an increasing sequence of smooth functions from $[-1,1]$ to $[0,1]$ with $\phi(x)=1$ for $x \in [-1+\frac 2n,1- \frac 2n]$ and $\phi(x)=0$ for $|x|\ge 1-\frac 1n$.  Then $(\phi_n(F_t))_{n \in \bbbn}$ is uniformly bounded and converges strongly to the identity for $n \to \infty$. Furthermore $\phi_n(F_t)\in l^1(\N)$ for $t \in [0,1]$. (Here is where we need the assumption that $F_t$ is the bounded transform of a selfadjoint operator with resolvents in $K(\N)$.) 

With $f:=e^{\pi i (\chi+1)}-1$
$$\tau (e^{- \pi i (\chi(F_t)+1)} \frac{d}{dt} f(F_t)) = \lim_{n \to \infty}\tau (e^{- \pi i (\chi(F_t)+1)}(\frac{d}{dt}f(F_t))\phi_n(F_t)) \ .$$

We use Lemma \ref{der} and get
\begin{eqnarray*}
\lefteqn{\tau (e^{- \pi i (\chi(F_t)+1)}(\frac{d}{dt}f(F_t))\phi_n(F_t))}\\
&=&\tau (e^{- \pi i (\chi(F_t)+1)}\int_{-\infty}^{\infty} \int_0^1 i \lambda \hat f(\lambda) e^{i u\lambda F_t}(\frac{d}{dt}F_t)e^{i (1-u) \lambda F_t}\phi_n(F_t)~ du~d\lambda)\\
&=&\tau ((\frac{d}{dt}F_t)f'(F_t)e^{- \pi i (\chi(F_t)+1)}\phi_n(F_t))\\
&=&i\pi ~\tau((\frac{d}{dt}F_t)\chi'(F_t)\phi_n(F_t))
\end{eqnarray*}
In the limit $n \to \infty$ we obtain that
$$\tau (e^{- \pi i (\chi(F_t)+1)} \frac{d}{dt} f(F_t))=i\pi ~\tau((\frac{d}{dt}F_t)\chi'(F_t)) \ .$$
Hence
$$\spfl((F_t)_{t \in [0,1]})-L = \frac{1}{2}  \int_0^1  \tau ((\frac{d}{dt}F_t)\chi'(F_t)) ~dt \ .$$
We evaluate $L$: 
We have that
\begin{eqnarray*}
\frac{1}{2\pi i}\int_0^1 \tau  (e^{-\pi i (G_t+1)}\frac{d}{dt}e^{\pi i (G_t+1)})~dt&=& \frac{1}{2}\int_0^1\tau(\frac{d}{dt}G_t)~dt\\
&=&-\frac 12\tau((2P_0-1)-\chi(F_0)) \ .
\end{eqnarray*}
Analogously
$$\frac{1}{2\pi i}\int_0^1 \tau  (e^{-\pi i (H_t+1)}\frac{d}{dt}e^{\pi i (H_t+1)})~dt=\frac 12\tau((2P_1-1)-\chi(F_1)) \ .$$
Thus
$$L= \frac 12 \tau((2P_1-1)- \chi(F_1))- \frac 12 \tau((2P_0-1)- \chi(F_0)) \ .$$
\end{proof}

In the following we derive a variation of the above theorem for paths of unbounded operators with common domain. Such a situation arises naturally in geometric examples (see \S \ref{L2ind}).    

For a closed operator $D$ on $H$ we denote by $H(D)$ the Hilbert space whose underlying vector space is $\dom D$ and whose scalar product is given by $$\langle x,y\rangle_D:=\langle x,y\rangle + \langle Dx,Dy\rangle \ .$$ 
By $B(H(D),H)$ we mean the space of bounded operators from $H(D)$ to $H$ endowed with the operator norm.

Let $(D_t)_{t \in [0,1]}$ be a path of operators in $ASF(\N)$ with common domain.

If $(t \mapsto D_t) \in C([0,1],B(H(D_0),H))$, then by \cite[Prop. 2.2]{l} the bounded transform $F_t=D_t(1+D_t^2)^{-1/2}$ depends continuously on $t$. However it is not clear to the author whether or not $(t \mapsto D_t) \in C^1([0,1],B(H(D_0),H))$ implies that $(t \mapsto F_t) \in C^1([0,1],\N)$, or at least $(t \mapsto F_tK) \in C^1([0,1],l^1(\N))$ for $K \in l^1(\N)$. Therefore we will modify the previous theorem.
 
Before doing so we discuss some technical points, which will be used in the following without further notice.

If $D,D' \in AS(\N)$ and $D' \in B(H(D),H)$, then $D'(1+D^2)^{-1/2} \in \N$.  This holds since $f_n(D')D'(1+D^2)^{-\frac 12} \in \N$ converges strongly to $D'(1+D^2)^{-1/2}$, if $f_n$ is defined as $f_n(x):=f(\frac xn)$ for some $f \in C_c(\bbbr)$ with $f|_{[-1,1]}=1$. 

In particular $\frac 1h(D_{t+h}-D_t)(D_t+\lambda)^{-1} \in \N$ for $h >0$ and $\lambda \notin \bbbr$ if $(t \mapsto D_t) \in C^1([0,1],B(H(D_0),H))$. It follows that
$$(\frac{d}{dt}D_t)(D_t+\lambda)^{-1} \in \N \ .$$
Furthermore $(t \mapsto (D_t-\lambda)^{-1}) \in C^1([0,1],B(H,H(D_0)))$ for $\lambda \notin \bbbr$ and 
$$\frac{d}{dt}(D_t-\lambda)^{-1}=-(D_t-\lambda)^{-1}(\frac{d}{dt}D_t)(D_t-\lambda)^{-1} \ .$$

We denote by $l^1(D_0)$ the Banach space $(D_0+i)^{-1}l^1(\N)$ with norm $\|T\|_{D_0}=\|(D_0+i)T\|_1$. Then the map $(D_0+i):l^1(D_0) \to l^1(\N)$ and its inverse are continuous.

In the statement of the following theorem $\chi$ is as defined before Theorem \ref{thintfor1}.

\begin{theorem}
Let $(D_t)_{t \in [0,1]}$ be a path of operators in $ASK(\N)$ with common domain such that $(t \mapsto D_t) \in C^1([0,1],B(H(D_0),H))$. Let $F_t:=D_t(1+D_t^2)^{-1/2}$.
Assume that $(t\mapsto \chi'(F_t)) \in C([0,1],l^1(D_0))$ and $(t\mapsto (\chi(F_t)^2-1)) \in C^1([0,1],l^1(\N))\cap C([0,1],l^1(D_0))$. Furthermore let $F_0,F_1$ be invertible.  

Then 
\begin{eqnarray*}
\spfl((D_t)_{t \in [0,1]})&=&\frac 12 \int_0^1  \tau ((\frac{d}{dt}F_t)\chi'(F_t))~ dt\\
&& +~ \frac 12 \tau(2P_1-1- \chi(F_1)) - \frac 12 \tau(2P_0-1-\chi(F_0))\ ,
\end{eqnarray*} 
where $P_i= 1_{\ge 0}(F_i)$.
\end{theorem}

\begin{proof}
First we adapt Lemma \ref{hol}, Prop. \ref{der} and Lemma \ref{holtrace} to the present situation.

Let $g \in C(\bbbr)$.
From $(t \mapsto F_t) \in C([0,1],\N)$ it follows that $(t \mapsto g(F_t)) \in C([0,1],\N)$.
Since $(D_t+i)(D_0+i)^{-1}$ as well as its inverse depend continuously on $t$ in $\N$, we have that $$(D_0+i)g(F_t)(D_0+i)^{-1}=(D_0+i)(D_t+i)^{-1}g(F_t)(D_t+i)(D_0+i)^{-1}$$ depends continuously on $t$ as well, hence for $K \in l^1(D_0)$ $$(t \mapsto g(F_t)K) \in C([0,1],l^1(D_0)) \ .$$ 

For $K\in l^1(D_0)$ 
$$F_tK=\frac{2}{\pi}\int_0^{\infty} (D_t^2+1+ \lambda^2)^{-1}D_t K ~d\lambda \ .$$
Using 
$$(D_t^2+1+ \lambda^2)^{-1}=(D_{t}+i\sqrt{1 + \lambda^2})^{-1}(D_{t}-i\sqrt{1 + \lambda^2})^{-1}$$ we get that  
\begin{eqnarray*}
\frac{d}{dt}F_tK &=&\frac{2}{\pi}\int_0^{\infty} (D_t^2+1+\lambda^2)^{-1}\frac{d}{dt}D_t K~ d\lambda \\
&&-~\frac{2}{\pi}\int_0^{\infty} (D_t^2+1 +\lambda^2)^{-1}(\frac{d}{dt}D_t)(D_{t}-i\sqrt{1 + \lambda^2})^{-1}D_t K ~ d\lambda \\
&&-~\frac{2}{\pi}\int_0^{\infty} (D_t+i\sqrt{1 +\lambda^2})^{-1}(\frac{d}{dt}D_t)(D_{t}^2+1+\lambda^2)^{-1}D_t K  ~ d\lambda \ .
\end{eqnarray*}
The integral converges in $l^1(\N)$.
Hence  for $K \in l^1(D_0)$
$$(t \mapsto F_tK) \in C^1([0,1],l^1(\N)) \ .$$ 

The proof of Lemma \ref{hol} modifies showing that for any holomorphic function $g$ defined in a neighbourhood of $[-1,1]$ and any $K \in l^1(D_0)$  
$$(t \mapsto g(F_t)K) \in C^1([0,1],l^1(\N)) \ .$$  

By adapting the proof of Prop. \ref{der} one obtains that for any $g \in C^2_c(\bbbr)$ and $K \in l^1(D_0)$
$$(t \mapsto g(F_t)K) \in C^1([0,1],l^1(\N))$$
and that the formula for the derivative from Prop. \ref{der} holds in this situation.
  
The function $(e^{\pi i (z+1)}-1)/(z^2-1)$ extends to an entire function $g$. Hence $$e^{\pi i (\chi(F_t)+1)}-1=g(\chi(F_t))(\chi(F_t)^2-1) \in C^1([0,1],l^1(\N)) \ .$$   
Using these modifications the proof of the theorem proceeds as the proof of the previous theorem.
\end{proof}

For the ordinary spectral flow a similar formula was derived in \cite{waspflfor}. 

In the following we derive generalizations of some of the integral formulas proven by Carey--Phillips \cite{cp1}\cite{cp2} for bounded perturbations of a fixed operator. We use their method of relating integral formulas for a path of unbounded operators to integral formulas for the path of bounded transforms.

\begin{lem}
\label{der2}
Let $f$ be a holomorphic function defined on $\{x+iy~|~(x,y)\in \bbbr \times ]-\ve,\ve[\} \subset \bbbc$ for some $\ve>0$ and assume that there are $C,c>0$ such that $$|f(\lambda)| \le C(1+|\lambda|)^{-c-1} \ .$$ Let $(D_t)_{t \in [0,1]}$ be a path of operators in $AS(\N)$ with common domain such that $(t \mapsto D_t) \in C^1([0,1],B(H(D_0),H))$. Let $\Gamma=(\bbbr-\frac{\ve}{2}i)\cup (\bbbr + \frac{\ve}{2}i)$, where we endow $\bbbr-\frac{\ve}{2}i$ with the standard orientation and $\bbbr+\frac{\ve}{2}i$ with the reversed orientation. Then $(t \mapsto f(D_t)) \in C^1([0,1],\N)$ and $$\frac{d}{dt}f(D_t)=-\frac{1}{2\pi i} \int_{\Gamma} f(\lambda)(D_t-\lambda)^{-1}(\frac{d}{dt}D_t)(D_t-\lambda)^{-1}~d\lambda \ .$$
\end{lem}

\begin{proof}
The assertion follows from Cauchy formula
$$f(D_t)=\frac{1}{2\pi i} \int_{\Gamma} f(\lambda)(D_t-\lambda)^{-1}~d\lambda \ .$$
\end{proof}

\begin{prop}
\label{thetasum}

Let $(D_t)_{t \in [0,1]}$ be a path of operators in $ASK(\N)$ with common domain such that $(t \mapsto D_t) \in C^1([0,1],B(H(D_0),H))$. We assume that there is $\delta \in (0,1/2)$ such that $e^{-\delta D_t^2} \in l^1(\N)$ for all $t \in [0,1]$ and $C>0$ such that $\|e^{-\delta D_t^2}\|_1<C$ for all $t \in [0,1]$. Furthermore let $D_0,D_1$ be invertible and set $P_i= 1_{\ge 0}(D_i)$.  

Let $F_t=D_t(1+D_t^2)^{-1/2}$.

Let $\chi_e(x)= \frac{1}{C}\int_0^x (1-y^2)^{-\frac 32}e^{(1-\frac{1}{1-y^2})}~dy$ with $C=\int_0^1 (1-y^2)^{-\frac 32}e^{(1-\frac{1}{1-y^2})}~dy=\frac{\sqrt{\pi}}{2}$.

Then 
\begin{eqnarray*}
\spfl((D_t)_{t \in [0,1]})&=&\frac{1}{\sqrt{\pi}} \int_0^1  \tau ((\frac{d}{dt}D_t)e^{-D_t^2})~ dt\\
&& +~ \frac 12 \tau(2P_1-1- \chi_e(F_1)) - \frac 12 \tau(2P_0-1-\chi_e(F_0))\ .
\end{eqnarray*} 
\end{prop}

\begin{proof} 
Let $(1-\ve)>2\delta$. First we show that 
$$(t\mapsto e^{-(1-\ve) D_t^2}) \in C^1([0,1],l^1(\N)) \ .$$ Lemma \ref{der2} implies that $(t \mapsto e^{-sD_t^2}) \in C^1([0,1],\N)$ for $s>0$. We have that
\begin{eqnarray*}
e^{-(1-\ve)D_{t+h}^2}-e^{- (1-\ve)D_t^2}&=&e^{- \delta D_{t+h}^2}(e^{-(1-\ve-\delta)D_{t+h}^2} - e^{-(1-\ve-\delta)D_t^2})\\
&& + (e^{-\delta D_{t+h}^2}-e^{-\delta D_t^2})e^{-(1-\ve-\delta) D_t^2} \ .
\end{eqnarray*}
Since $\|e^{-\delta D_t^2}\|_1$ is uniformly bounded, the limit of the right hand side is zero for $h \to 0$, thus 
$$(t \mapsto e^{-(1-\ve)D_t^2}) \in C([0,1],l^1(\N)) \ .$$ Using this it also follows that the limit for $h \to 0$ of the right hand side divided by $h$ exists. Thus 
$$(t \mapsto e^{-(1-\ve)D_t^2}) \in C^1([0,1],l^1(\N)) \ .$$

Now let $\chi \in C^2_c(\bbbr)$ be a function that equals $\chi_e$ in a neighbourhood of $[-1,1]$.

We have that $1-F_t^2=(1+D_t^2)^{-1/2}$, hence 
$$(1-F_t^2)^{-3/2}e^{(1-\frac{1}{1-F_t^2})}=(1+D_t^2)^{3/2}e^{-D_t^2} \ .$$

From $(t \mapsto (D_0+i)(1+D_t^2)^{3/2}e^{-\ve D_t^2}) \in C([0,1],\N)$ it follows that $(t \mapsto \chi'(F_t)) \in C([0,1],l^1(D_0))$.

The function $\psi(x):=(1-\chi(x)^2)e^{-(1-\frac{2\ve}{3})(1-\frac{1}{1-x^2})}$, defined on $(-1,1)$, extends to an element in $C_c^2(\bbbr)$. Thus, by the modification of Prop. \ref{der} formulated in the proof of the previous theorem, for $K\in l^1(D_0)$
$$(t \mapsto \psi(F_t)K) \in C([0,1],l^1(D_0)) \cap C^1([0,1],l^1(\N)) \ .$$ From $$(1-\chi(F_t)^2)=\psi(F_t)e^{-(1-\frac{2\ve}{3})D_t^2}$$ it follows that $$(t \mapsto (1-\chi(F_t)^2)) \in C([0,1],l^1(D_0)) \cap C^1([0,1],l^1(\N)) \ .$$

Hence we can apply the previous theorem and get
\begin{eqnarray*}
\spfl((D_t)_{t \in [0,1]})&=&\frac{1}{2C} \int_0^1  \tau ((\frac{d}{dt}F_t)(1+D_t^2)^{3/2}e^{-D_t^2})~ dt\\
&& +~ \frac 12 \tau(2P_1-1- \chi(F_1)) - \frac 12 \tau(2P_0-1-\chi(F_0))\ .
\end{eqnarray*}

The equality 
$$\tau ((\frac{d}{dt}F_t)(1+D_t^2)^{3/2}e^{-D_t^2})=\tau ((\frac{d}{dt}D_t)e^{-D_t^2})$$ 
can be shown in a straightforward way as follows: One expresses the right hand side of $$(\frac{d}{dt}F_t)(1+D_t^2)^{3/2}e^{-D_t^2}=\frac{d}{dt}(F_t(1+D_t^2)^{3/2}e^{-D_t^2})-F_t\frac{d}{dt}((1+D_t^2)^{3/2}e^{-D_t^2})$$ in terms of $D_t$ and $\frac{d}{dt}D_t$ by applying Lemma \ref{der2} and then one uses the cyclicity of the trace.
\end{proof}

In the following lemma we give a different expression for the contributions from the endpoints, which will be needed in the proof of the subsequent proposition.

\begin{lem} 
\label{eta}
Let $D \in ASK(\N)$ be invertible and let $\chi \in C^1(\bbbr)$ be a normalizing function for $D$ such $\chi(D)^2-1 \in l^1(\N)$ and such that $(t \mapsto \chi'(\sqrt t D))$ is an integrable function in $C([0,\infty),l^1(D))$.  Then with $P=1_{\ge 0}(D)$
 $$\tau(2P-1- \chi(D))=\frac 12\int_1^{\infty}t^{-1/2}\tau(D\chi'(\sqrt t D))~dt \ .$$
\end{lem}

\begin{proof}
The right hand side of the equation is well-defined since $(2P-1)-\chi(D)=((2P-1)+\chi(D))^{-1}(1-\chi(D)^2) \in l^1(\N)$.

For $x \neq 0$
$$\sign(x)- \chi(x)=\frac 12\int_1^{\infty}t^{-1/2} x\chi'(\sqrt t x)~dt \ .$$

This proves that $2P-1- \chi(D)=\frac 12\int_1^{\infty}t^{-1/2} D\chi'(\sqrt t D)~dt$. The assertion follows since the integral converges in $l^1(\N)$.
\end{proof}

For invertible $D$ such that $e^{-(1-\ve)D^2} \in l^1(\N)$ for some $\ve>0$ the truncated $\eta$-invariant (see \cite[\S 8]{cp2}) is defined by
$$\eta_1(D)=\frac{1}{\sqrt \pi}\int_1^{\infty} t^{-1/2} \tau(D e^{-tD^2})~dt \ .$$

The previous lemma applied to $\chi(x)=\chi_e(x(1+x^2)^{-1/2})$ shows that
$$\tau(2P-1-\chi_e(F_D))=\eta_1(D) \ .$$

We define for $p \ge 1$ the ideal
$$l^p(\N)=\{A \in \N~| ~|A|^p \in l^1(\N)\}$$ with norm $\|A\|_p=\tau(|A|^p)^{1/p}+\|A\|$. For $A \in l^p(\N)$ and $S,T \in \N$ 
$$\|SAT\|_p \le \|S\|~ \|A\|_p~ \|T\| \ .$$
We refer to \cite{d2} for proofs of these facts. By \cite[Prop. A.1]{cp0} the space $l^p(\N)$ is a Banach space. 

\begin{prop}
\label{psum}
Let $(D_t)_{t \in [0,1]}$ be a path of operators in $ASK(\N)$ with common domain such that $(t \mapsto D_t) \in C^1([0,1], B(H(D_0),H))$. Furthermore assume that $D_0,D_1$ are invertible and set $P_i=1_{\ge 0}(D)$. We also assume that there is $p \in [1,\infty)$ such that $(1+D_0^2)^{-p/2} \in l^1(\N)$.

We write $F_t:=D_t(D+D_t^2)^{-1/2}$.

Define $\chi_p(x)=\frac{1}{C_p} \int_0^x (1-y^2)^{(p-2)/2}~dy$ with $C_p=\int_0^1(1-y^2)^{(p-2)/2}~dy$.

Then 
\begin{eqnarray*}
\spfl((D_t)_{t \in [0,1]}) &=& \frac{1}{2C_p}\int_0^1 \tau((\frac{d}{dt} D_t)(1+D_t^2)^{-(p+1)/2})~dt \\
&& + \frac 12 \tau(2P_1-1-\chi_p(F_1))-\frac 12 \tau(2P_0-1-\chi_p(F_0)) \ .
\end{eqnarray*}
\end{prop}

\begin{proof}
We use the method of \cite[\S 9]{cp2}. Some of the following arguments can be found there in more detail.

First we prove that $\tau(2P_i-1-\chi_p(F_i))$ is well-defined:
Using L'Hospital rule it can be checked that $(\chi_p(x)^2-1)(1-x^2)^{-p/2}$ extends to a bounded continuous function on $[-1,1]$, hence $(\chi_p(F_i)^2-1)(1-F_i^2)^{-p/2} \in \N$. From $(1-F_i^2)^{p/2}=(1+D_i^2)^{-p/2} \in l^1(\N)$ it follows that $(\chi_p(F_i)^2-1) \in l^1(\N)$.

Since $(D_t \pm i)^{-1}(D_0+ i)$ depends continuously on $t$ in $\N$ and $(D_0 + i)^{-1} \in l^p(\N)$ we have that $(D_t \pm i)^{-1}$ depends continuously on $t$ in $l^p(\N)$ and hence that $\tau((1+D_t^2)^{-p/2})$ depends continuously on $t$. By $$e^{-sD_t^2}=e^{-sD_t^2}(1+D_t^2)^{p/2} (1+D_t^2)^{-p/2} \ ,$$ for $s>0$ fixed $e^{-sD_t^2}$ is uniformly bounded in $t$ as an element in $l^1(\N)$. Thus for $s>0$ by Cor. \ref{thetasum}
$$\spfl((D_t)_{t \in [0,1]})=\spfl((\sqrt s D_t)_{t \in [0,1]})= \sqrt{\frac{s}{\pi}} \int_0^1  \tau ((\frac{d}{dt}D_t)e^{-s D_t^2})~ dt \ .$$
Using $$\int_0^{\infty} s^{(p-2)/2}e^{-s}~ds=\Gamma(\tfrac p2) \ ,$$
where $\Gamma$ denotes the Gamma-function, we get that
\begin{eqnarray*}
\spfl((D_t)_{t \in [0,1]})&=&\frac{1}{\Gamma(\frac p2)} \int_0^{\infty} s^{(p-1)/2}~ \frac{1}{\sqrt{\pi}}\int_0^1  \tau ((\frac{d}{dt}D_t)e^{-s(1+ D_t^2)})~ dt\dots \\
&&\dots + \frac 12 s^{(p-2)/2}e^{-s}(\eta_1(\sqrt s D_1)- \eta_1(\sqrt s D_0))~ds \ .
\end{eqnarray*}

We study the right hand side:

For $q \ge \frac 12$
$$\int_0^{\infty} s^q (1+D_t^2)^{1/2} e^{-s(1+D_t^2)}~ds=\Gamma(q+1)(1+D_t^2)^{-q-1/2} \ ,$$
where the integral converges in the operator norm. Let $q\ge\frac{(p-1)}{2}$. Since $s^q (1+D_t^2)^{1/2}e^{-s(1+D_t^2)}$ is positive and in $l^1(\N)$ for $s \in (0,\infty)$ and since the right hand side is in $l^1(\N)$, the integral converges also in $l^1(\N)$. 
Hence 
$$(s \mapsto s^q e^{-s(1+D_t^2)}) \in L^1(\bbbr^+,l^1(D_0)) \ .$$
Using Fubini theorem we can also conclude that
$$((s,t) \mapsto s^q \tau((1+D_t^2)^{1/2} e^{-s(1+D_t^2)})) \in L^1(\bbbr^+ \times [0,1]) \ .$$
By
$$|\tau((\frac{d}{dt}D_t)e^{-s(1+D_t^2)})|\le \| (\frac{d}{dt} D_t)(1+D_t^2)^{-1/2}\|~ \tau((1+D_t^2)^{1/2}e^{-s(1+D_t^2)})$$
this implies that $$(s,t) \mapsto s^q\tau((\frac{d}{dt}D_t)e^{-s(1+D_t^2)}) \in L^1(\bbbr^+ \times [0,1]) \ .$$
Thus
\begin{eqnarray*}
\lefteqn{\frac{1}{\Gamma(\frac p2)\sqrt{\pi}}\int_0^{\infty}\int_0^1 s^{(p-1)/2}\tau ((\frac{d}{dt}D_t)e^{-s(1+ D_t^2)})~ dtds}\\
&=&\frac{1}{\Gamma(\frac p2)\sqrt{\pi}} \int_0^1 \tau \bigl((\frac{d}{dt}D_t)\int_0^{\infty}s^{(p-1)/2} e^{-s(1+ D_t^2)})~ds\bigr)~dt\\
&=&\frac{\Gamma(\frac{p+1}{2})}{\Gamma(\frac p2)\sqrt{\pi}}\int_0^1 \tau((\frac{d}{dt}D_t)(1+D_t^2)^{-(p+1)/2}) dt \ .
\end{eqnarray*}

Now consider the contributions of the endpoints. Let $i=0,1$.

Since $D_i$ is invertible, 
$$((s,t) \mapsto s^{(p-1)/2}t^{-1/2} D_i e^{-s(1+ tD_i^2)}) \in L^1(\bbbr^+ \times [1,\infty),l^1(\N)) \ .$$
We evaluate
\begin{eqnarray*}
\lefteqn{\frac{1}{2\Gamma(\frac p2)}\int_0^{\infty}s^{(p-2)/2}e^{-s}\eta_1(\sqrt s D_i)~ds}\\
&=& \frac{1}{2\Gamma(\frac p2)\sqrt{\pi}}\int_0^{\infty}\int_1^{\infty} s^{(p-1)/2}t^{-1/2}\tau(D_i e^{-s(1+ tD_i^2)})~dtds\\
&=&\frac{\Gamma(\frac{p+1}{2})}{2\Gamma(\frac p2)\sqrt{\pi}}\int_1^{\infty} t^{-1/2}\tau(D_i (1+t D_i^2)^{-(p+1)/2})~dt \ .
\end{eqnarray*}
Now we apply the previous lemma to the normalizing function $\chi(x):=\chi_p(x(1+x^2)^{-1/2})$ of $D$. Since $\chi'(x)=\frac{1}{C_p} (1+x^2)^{-(p+1)/2}$, the last line equals
$$\frac{C_p\Gamma(\frac{p+1}{2})}{2\Gamma(\frac p2)\sqrt{\pi}}\int_1^{\infty} t^{-1/2}\tau(D_i\chi'(\sqrt t D_i))~dt=\frac{C_p\Gamma(\frac{p+1}{2})}{\Gamma(\frac p2)\sqrt{\pi}} \tau(2P_i-1-\chi_p(F_i)) \ .$$ 

The assertion follows from
$$\frac{\Gamma(\frac p2)\sqrt{\pi}}{\Gamma(\frac{p+1}{2})}=\int_0^1x^{(p-2)/2}(1-x)^{-1/2}~dx=2C_p \ ,$$
where the first equality follows from the properties of the Beta-function and the last equality is obtained from the change of variable $x=(1-y^2)$.
\end{proof}

\section{$L^2$-index theory}
\label{L2ind}

In this section we will show that the integral formulas of Prop. \ref{thetasum} and Prop. \ref{psum} apply to paths of invariant symmetric elliptic differential operators on a Galois covering. The proofs rely on the theory of regular operators on Hilbert $C^*$-modules for which we refer to \cite{la} and the theory of pseudodifferential operators over $C^*$-algebras, which was developed in \cite{mf}. Most of the material assembled in this section is well-known at least to experts.

\begin{lem}
Let $\A,\B$ be $C^*$-algebras and $H_1$ resp. $H_2$ a Hilbert $\A$-module resp. Hilbert $\B$-module. Let $\rho:\A \to B(H_2)$ be a $C^*$-homomorphism and let $D_1$ be a regular selfadjoint operator on $H_1$.

\begin{enumerate}
\item There is an induced regular selfadjoint operator $D$ on $H_1 \ten_{\rho} H_2$ such that the homomorphism $\rho_*:B(H_1) \to B(H_1 \ten_{\rho} H_2), ~ T \mapsto T \ten 1$ maps $f(D_1)$ to $f(D)$ for all $f\in C(\bbbr)$ for which $\lim_{x \to \infty}f(x)$ and $\lim_{x \to -\infty}f(x)$ exist.
 
\item Let $S$ be a core of $D_1$. Then the span of the set $\{x \ten y \in H_1 \ten_{\rho} H_2~|~x \in S,~y \in H_2\}$ is a core of $D$, and $D$ acts on it by $D(x \ten y)=(D_1x) \ten y$.    

\item Denote by $H(D_1)$ the Hilbert $\A$-module whose underlying $\A$-module is $\dom D_1$ and whose $\A$-valued scalar product is given by $$\langle x,y \rangle_{D_1}:=\langle x,y \rangle +\langle D_1x,D_1y \rangle \ .$$ Define analogously $H(D)$.

The map $H(D_1) \ten_{\rho} H_2 \to H(D)$ that maps $x \ten y,~x \in S,y\in H_2$ to $x \ten y \in \dom D$, is an isomorphism.
\end{enumerate} 
\end{lem}

\begin{proof}
As usual let $F_{D_1}=D_1(1+D_1^2)^{-\frac 12}$. Since $\Ran(1-F_{D_1}^2)=\Ran ((1+D_1^2)^{-1})$ is dense in $H_1$, the range of $\rho_*(1- F_{D_1}^2)$ is dense in $H_1 \ten_{\rho} H_2$. By \cite[Theorem 10.4]{la} it follows that $\rho_*(F_{D_1})$ is the bounded transform of a regular selfadjoint operator $D$ on $H_1 \ten_{\rho} H_2$. 
This shows (1).

Since $D_1$ and $D$ are regular, the operators $(D_1 +i) \ten 1:H(D_1)\ten_{\rho} H_2 \to H_1 \ten_{\rho} H_2$ and $D+i:H(D) \to H_1 \ten_{\rho} H_2$ are isomorphisms. The isomorphism $U:=(D+i)^{-1}(D_1 \ten 1 +i):H(D_1) \ten_{\rho} H_2 \to H(D)$ is given by $v \ten w \mapsto v\ten w$. Hence (3) follows. 

The operator $D_1 \ten 1 + i$ is determined by its action on the set $\{x \ten y~|~x \in S,~y \in H_2\}$. By $(D+i)=(D_1\ten 1 +i)U^{-1}$ and since $U$ preserves this set, the operator $D$ is also determined on this set.
\end{proof}

Let $M$ be a closed Riemannian manifold and let $p:\tilde M \to M$ be a Galois covering of $M$ with deck transformation group $\Gamma$. The manifold $\tilde M$ inherits an invariant Riemannian metric from $M$. Let $E$ be a hermitian vector bundle on $M$ and let $\tilde E = p^*E$ be endowed with the invariant hermitian structure induced by $E$. We identify $E$ with its dual.

In the following we recall the construction of the correspondence of invariant differential operators on $\C_c(\tilde M,\tilde E)$ and differential operators on $\C(M,E\ten \Pj)$, where $\Pj$ is the Mishchenko-Fomenko vector bundle.
 More details on the following discussion can be found, for example, in \cite[\S 3]{waindfor} and \cite[\S E.3]{ps}.

Let $C_r^*\Gamma$ be the reduced group $C^*$-algebra. Recall that $\Pj$ is the $C_r^*\Gamma$-vector bundle $\tilde M \times_{\Gamma} C_r^*\Gamma$ and let $\Pj_{alg}=\tilde M \times_{\Gamma} \bbbc\Gamma$. The right $\Gamma$-action $R$ on $\tilde E$ induces a right $\Gamma$-action $$\Gamma \times \C_c(\tilde M,\tilde E) \to \C_c(\tilde M,\tilde E),~(g,s) \mapsto R_{g^{-1}}^*s$$ with $(R_{g^{-1}}^*s)(x):=R_{g}(s(xg^{-1}))$. 

Furthermore there are a left and a right $\Gamma$-action $$\Gamma \times \C(\tilde M, \tilde E \ten \bbbc \Gamma) \to \C(\tilde M, \tilde E \ten \bbbc \Gamma) \ :$$ 
The left $\Gamma$-action is given by $$(h,\sum_{g \in \Gamma}s_g g) \mapsto \sum_{g \in \Gamma}(R_h^*s_g) hg$$
and the right $\Gamma$-action by 
$$(h,\sum_{g \in \Gamma}s_g g) \mapsto \sum_{g \in \Gamma}s_g gh \ .$$
The invariant subspace $\C(\tilde M, \tilde E \ten \bbbc \Gamma)^{\Gamma}$ with respect to the left $\Gamma$-action  is isomorphic to $\C(M,E\ten \Pj_{alg})$. The isomorphism is $\Gamma$-equivariant with respect to the right $\Gamma$-action.
Moreover the map 
$$\C_c(\tilde M,\tilde E) \to  \C(\tilde M, \tilde E \ten \bbbc \Gamma)^{\Gamma},~ s\mapsto \sum_{g \in \Gamma}(R_g^*s) g$$  is a $\Gamma$-equivariant isomorphism. Hence there is a $\Gamma$-equivariant isomorphism $$\C_c(\tilde M,\tilde E) \cong \C(M,E\ten \Pj_{alg}) \ .$$ 

Let $D:\C_c(\tilde M,\tilde E) \to \C_c(\tilde M,\tilde E)$ be an invariant symmetric elliptic differential operator. We denote its closure as an unbounded operator on $L^2(\tilde M,\tilde E)$ by $D$ as well.  Via the previous isomorphism the operator $D$ induces an elliptic symmetric differential operator $\D:\C(M,E\ten \Pj_{alg}) \to \C(M,E\ten \Pj_{alg})$. Its closure on the Hilbert $C^*_r\Gamma$-module $L^2(M,E \ten \Pj)$, denoted by $\D$ as well, is a regular selfadjoint operator. We define $\rho:C^*_r\Gamma \to B(l^2(\Gamma))$ to be the inclusion. By applying the previous lemma we get a selfadjoint operator $D'=\D \ten 1$ on the Hilbert space $L^2(M,E \ten \Pj)\ten_{\rho} l^2(\Gamma)$.

The $\Gamma$-equivariant isomorphism $$\C(M,E\ten \Pj_{alg})\ten_{\bbbc \Gamma}\bbbc\Gamma \cong \C(M,E\ten \Pj_{alg})\cong \C_c(\tilde M,\tilde E)$$ extends to a $\Gamma$-equivariant isometry
$$L^2(M,E \ten \Pj)\ten_{\rho} l^2(\Gamma) \cong L^2(\tilde M,\tilde E) \ ,$$
which intertwines $D$ with $D'$.

By the previous lemma, one can use information on $\D$ in order to gain knowledge about $D$ (see also \cite[Prop. E.6]{ps}).

\begin{prop}
\label{diffell}
Let $(D_t)_{t\in [0,1]}$ be a path of invariant selfadjoint elliptic differential operators of order $k>0$ on $L^2(\tilde M,\tilde E)$ with coefficients that are of class $C^1$ in $t$. We have that $\dom D_t=\dom D_0$ for all $t \in [0,1]$ and $(t \mapsto D_t) \in C^1([0,1],B(H(D_0),H))$ for $H=L^2(\tilde M,\tilde E)$.
\end{prop}

\begin{proof}
Let $(\D_t)_{t \in [0,1]}$ be the corresponding family of elliptic differential operators on $L^2(M,E \ten \Pj)$. The coefficients of $\D_t$ are of class $C^1$ in $t$ as well. Furthermore $H(\D_0)$ equals the $k$-th Sobolev space $H^k(M,E \ten \Pj)$. Hence 
$$(t \mapsto \D_t) \in C^1([0,1],B(H(\D_0),L^2(M,E \ten \Pj))) \ .$$ By the first statement of the previous lemma this implies that
$$(t \mapsto D_t) \in C^1([0,1],B(H(\D_0) \ten_{\rho} l^2(\Gamma),L^2(\tilde M,\tilde E))) \ .$$ Since $H(\D_0) \ten_{\rho} l^2(\Gamma)$ can be identified with $H(D_0)$ by the third statement of the previous lemma, this implies the assertion.  
\end{proof}

In the following let $(D_t)_{t \in [0,1]}$ be as in the statement of the proposition.

We denote by $\N$ the von Neumann algebra of $\Gamma$-equivariant operators in $B(L^2(\tilde M,\tilde E))$.
Since $D_t$ commutes with the $\Gamma$-action, it is affiliated to $\N$.

Recall that the vector bundle $E \boxtimes E$ on $M\times M$ has fiber $(E \boxtimes E)_{(x,y)}=E_x \ten E_y$.

We consider $\Pj$ as isometrically embedded into $M \times (C_r^*\Gamma)^n$.

Let $\tr_{\Gamma}:E_x \ten E_x \ten M_n(C_r^*\Gamma) \to \bbbc$ be the trace induced by the trace $$C_r^*\Gamma \to \bbbc, ~\sum_{g \in\Gamma} s_g g\mapsto s_e $$ and the standard trace on $E_x \ten E_x\ten M_n(\bbbc)$. 

Let $$K:C(M \times M,(E \boxtimes E) \ten M_n(C_r^*\Gamma)) \to B(L^2(M,E \ten (C_r^*\Gamma)^n))$$ be the map that associates to an integral kernel the corresponding integral operator. 

The homomorphism $\rho$ induces $$\rho_*:B(L^2(M,E \ten (C_r^*\Gamma)^n)) \to \M:=B(L^2(M,E) \ten l^2(\Gamma)^n)^{\Gamma} \ .$$ 
We have that $\N \subset \M$.

The trace $\tau$ is defined by
$$\tau(\rho_*K(k))=\int_M \tr_{\Gamma} k(x,x)~dx \ ,$$
where $K(k)=K(k_1)K(k_2)$ with $k_1,k_2 \in C(M \times M,(E \boxtimes E) \ten M_n(C_r^*\Gamma))$.

It extends to a semifinite normal faithful trace on $\M$, which restricts to a semifinite normal faithful trace on $\N$. (This well-known fact can be shown by using the concept of Hilbert algebras, see \cite[\S I.5]{d} and \cite[\S I.6.2]{d}.)

There is an induced continuous map
$$\rho_* \circ K:C(M \times M,(E \boxtimes E) \ten M_n(C_r^*\Gamma)) \to l^2(\M) \ .$$

Furthermore the map $\rho_*$ maps ordinary trace class operators on $L^2(M,E^n)$ tensored with the identity on $M_n(C_r^*\Gamma)$ to $l^1(\M)$.

\begin{prop}
For $pk>\dim M$
$$(1+D_0^2)^{-p/2} \in l^1(\N) \ .$$
\end{prop}

\begin{proof}
Let $P:M \times (C_r^*\Gamma)^n \to \Pj$ be the orthogonal projection. Choose an elliptic symmetric pseudodifferential operator $R$ of order $k$ on $\C(M,E)$. By tensoring with the identity on $(C_r^*\Gamma)^n$ it defines a continuous operator from $L^2(M,E \ten (C_r^*\Gamma)^n)$ to $H^{kp}(M,E \ten (C_r^*\Gamma)^n)$. Set $$T:=\left(P(1+\D_0^2)^{-p/2}P+(1-P)(1+R^2)^{-p/2}(1-P)\right)(1+R^2)^{p/2} \ .$$ 
Since $T$ is a pseudodifferential operator of order zero, the operator $\rho_*(T)$ is an element in $\M$. Furthermore $(1+R^2)^{-p/2}$ is a trace class operator on $L^2(M,E)$, hence its image under $\rho_*$ is in $l^1(\M)$. It follows that the image of $$PT(1+R^2)^{-p/2}P=(1+\D_0^2)^{-p/2}$$ 
under $\rho_*$ is in $l^1(\N)$.
\end{proof}

\end{document}